\documentclass[12pt]{article}

\usepackage[T2A]{fontenc}
\usepackage[utf8]{inputenc}

\usepackage{amsfonts}
\usepackage{amssymb}
\usepackage{amsmath}
\usepackage{amsthm}

\def \le {\leqslant}
\def \ge {\geqslant}

\theoremstyle{plain}

\begin{document}
\begin{Large}
\centerline{
Einige Bemerkungen über  inhomogene diophantische Approximationen}
\end{Large}
\vskip+0.5cm
\begin{large}
\centerline{ Nikolay Moshchevitin}
\end{large}

\vskip+1cm

Es ist bekannt, dass
für jedes reelle irrationale  $\theta$ und  für jedes reelle  $\alpha$ 
unendlich viele Gitterpunkte $ (x,y)\in \mathbb{Z}^2, x\ge  1$ mit
$$
|\theta x - \alpha - y| \le \frac{1}{x}
$$
existieren.
Es gibt viele klassische mehrdimensionale Verallgemeinerungen dieser Aussage (siehe, zum Beispiel, 
 Kapitel V aus dem Buch von Cassels \cite{C}).
Das Ziel dieser Arbeit ist es, einige der  einfachsten Aussagen über inhomogene  diophantische Approximationen zu diskutieren. Außerdem denken wir, dass,    einige von diesen Ergebnissen (Sätze 3 und 7) möglicherweise nie dokumentiert wurden, obwohl die Ergebnisse sehr einfach sind. 
Zumindest sind uns keine Lehrbücher oder Artikel bekannt, in denen diese Ergebnisse explizit erwähnt werden.

Bezüglich moderner Ergebnisse im Zusammenhang mit inhomogenen Approximationen können wir, zum Beispiel, die Artikel 
  \cite{BL,Ein,F, Mon,MoMMJ}
erwähnen. Man kann
sehr interessante Ergebnisse, die sich mit den gewichteten Approximationen befassen,  in den Artikeln \cite{Ch,Ge} finden.
 
 \vskip+0.5cm

{\bf 1. Notationen.}
 \vskip+0.5cm

Seien $n$ und $ m$ natürliche Zahlen, $ d=m+n$.
Sei
$$\Theta
=\left(
\begin{array}{ccc}
\theta_{1,1}&\cdots&\theta_{1,m}\cr
\theta_{2,1}&\cdots&\theta_{2,m}\cr
\cdots &\cdots &\cdots \cr
\theta_{n,1}&\cdots&\theta_{n,m}
\end{array}
\right)\,\,\,
\,\,\,
\text{eine reele}\,\,
m\times n \,\text{ Matrix und sei}
$$
$$
\overline{\Theta} = 
 \left(
\begin{array}{ccccccc}
\theta_{1,1}&\cdots&\theta_{1,m}& -1 & 0&\cdots &0\cr
\theta_{2,1}&\cdots&\theta_{2,m}& 0& -1& \cdots &0\cr
\cdots &\cdots &\cdots &\cdots&\cdots&\cdots&\cdots\cr
\theta_{n,1}&\cdots&\theta_{n,m}&0&0&\cdots&-1
\end{array}
\right)\,\,\,
\,\,\,
\text{eine reele}\,\,
d\times n \,\text{ Matrix}.
$$
Für eine Matrix $A$ 
wird die transponierte Matrix mit $A^\top$
  bezeichnet. Wir betrachten die Vektoren
  $$
  \pmb{\theta}_j =
\left( \theta_{1,j},\dots,
\theta_{i,j},\dots,
\theta_{n,j}  \right)^\top\in \mathbb{R}^n,\,\,\,  j =1,...,m;
\,\,\,\,\,\,\,\,\,\,
\pmb{\alpha} =
\left(\alpha_1,
,\dots,
\alpha_n \right)^\top \in \mathbb{R}^n.
$$
und
$$
\pmb{x} =
\left(
x_{1},\cdots
,
x_{m} 
\right)^\top \in \mathbb{R}^m,
\,\,\,\,\,\,\,\,\,\,\,\,
\pmb{y} =
\left(
y_{1},\cdots,
y_{n} 
\right)^\top\in \mathbb{R}^n,
$$
$$
\pmb{z} =
\left(
z_1,\dots,
z_{d} 
\right)^\top
=
\left(
\pmb{x},
\pmb{y}
\right)^\top
=
\left(
x_{1},\dots,
x_{m},
y_{1},
\dots
,
y_{n} 
\right)^\top
\in \mathbb{R}^d
$$
Für $ \pmb{w} =(w_1,\dots, w_l)\in \mathbb{R}^l$ 
betrachten wir die sup-Norm $ |\pmb{w}| = \max_{1\le j \le l} |w_j|$ und den
Abstand 
$$
||\pmb{w}|| = \max_{1\le j \le l} ||w_j||,\,\,\,\,\,  ||w_j|| = \min_{a\in\mathbb{Z}} |w_j - a|
$$
zum nächsten Gitterpunkt. Durch
$$
\psi_{\Theta^\top} (t) =
\min_{\pmb{y}\in \mathbb{Z}^n: 0<|\pmb{y}|\le t}\,\,\,
||\Theta^\top \pmb{y}||
$$
bezeichnen wir
die 
Funktion des 
Irrationalitätsmaßes für die transponierte Matrix $\Theta^\top$.

 \vskip+0.5cm

{\bf 2. Über den Approximationssatz von Kronecker.}

 \vskip+0.5cm

Die folgende Aussage wurde von Khintchine \cite{H} bewiesen. Sie ist im Buch von Cassels  (siehe \cite{C}, Kapitel V, Theorem  XVII, B) angegeben.

   \vskip+0.3cm
{\bf Satz 1}.
{\it Nehmen wir an, dass  die Ungleichung 
\begin{equation}\label{2}
|| \alpha_1u_1+...+\alpha_n u_n|| \le  K \,\max
\left(
X \, ||\Theta^\top \pmb{u}||,  C\, |\pmb{u}|
\right),\,\,\,\,\, K= \frac{2^{n-1}}{(d!)^2}
\end{equation}
für alle $\pmb{u}\in\mathbb{Z}^n$
gilt.
Dann
gibt es einen Gitterpunkt $ \pmb{x}\in \mathbb{Z}^m$
mit
\begin{equation}\label{1}
||\Theta \pmb{x}-\pmb{\alpha}|| \le C ,\,\,\,\,\, |\pmb{x}| \le X.
\end{equation}
 
}

 \vskip+0.3cm
 Die folgende Aussage ist ein Klassisches Ergebnis.
  \vskip+0.3cm

{\bf Satz 2} (Approximationssatz von Kronecker).

{\it
Nehmen wir an, dass
$\pmb{y}\Theta \not\in \mathbb{Z}^m$ 
 für jedes $\pmb{y}\in \mathbb{Z}^n\setminus\{\pmb{0}\}$  gilt.
Dann gibt es für jedes $\varepsilon>0 $  und jedes $\pmb{\alpha} \in \mathbb{R}^n$ einen Gitterpunkt 
$ \pmb{z}  =
 \left(\begin{array}{c}\pmb{x}\cr \pmb{y}\end{array}\right)\in \mathbb{Z}^{m+n}$ 
 mit
$$
|\overline{\Theta}\pmb{z} - \pmb{\alpha} |=
|| {\Theta}\pmb{x} - \pmb{\alpha} ||
 < \varepsilon.
$$
}
 
   \vskip+0.3cm
 
 {\bf Bemerkung 1} (siehe \cite{C},  Kapitel 5 \S 8).
 Wir wollen ein paar Worte über den Beweis von Satz 2  sagen.
  Satz 2 folgt sofort aus Satz 1, weil     die Ungleichung
 \begin{equation}\label{000}
 \psi_{\Theta^\top} (t) >0
 \end{equation}
  für alle $  t\ge 1$
  gilt. 
  Sei $ C=\varepsilon^{-1}$
 Dann gibt es für jedes $C$ ein $X$ mit  
 $\psi_{\Theta^\top} \left(\frac{1}{2KC}\right) \ge \frac{1}{2KX}$.
 Oder, mit anderen Worten, 
 man hat 
$
||\Theta^\top \pmb{u}|| \ge\frac{1}{2KX}$
 für jedes $\pmb{u}\in \mathbb{Z}^n$ mit
$0\neq |\pmb{u}|\le \frac{1}{2KC}$.
Das bedeutet, dass 
die Ungleichung (\ref{2}) für jedes $\pmb{u}\in \mathbb{Z}^n$ gilt, sodass  es ein $\pmb{x}\in \mathbb{Z}^m$
mit  (\ref{1})  gibt.

   \vskip+0.3cm
 
 Seltsamerweise sind wir in der Literatur nie auf das folgende Ergebnis gestoßen.
    \vskip+0.3cm
 
{\bf Satz 3}.  
{\it
Wir betrachten den Fall $m \ge 2, n=1$.
Sei
$
\Theta = (\theta_1,...,\theta_m).
$
Nehmen wir an,
es gibt $i$ und $j$  mit $1\le i <j\le n$, sodass
 $ 1 , \theta_i,\theta_j$   linear unabhängig  sind über $\mathbb{Q}$.
Dann gibt es
 für jedes $\varepsilon>0 $  und jedes ${\alpha} \in \mathbb{R}$  ein $ \pmb{z}  =
 \left(\begin{array}{c}\pmb{x}\cr y\end{array}\right)\in \mathbb{Z}^{m+1}$ mit
$$
|\overline{\Theta}\pmb{z} - \pmb{\alpha} |=
||\theta_1 x_1+...+\theta_mx_m - \alpha || < \frac{\varepsilon}{|\pmb{x}|}.
$$
}

Beweis von Satz 3.  

V. Jarn\'{\i}k  \cite{J}  hat  bewiesen, dass  die asymptotische Gleichung
\begin{equation}\label{J}
\limsup_{t\to \infty} t \cdot \psi_{\Theta^\top} (t) = +\infty
\end{equation}
gilt (siehe Satz 9 aus \cite{J}, dieses Ergebnis wird auch in  \cite{M} diskutiert). Also gibt es
für jedes $\varepsilon >0$
 eine Folge $ t_\nu \to +\infty$ mit 
\begin{equation}\label{00}
t_\nu \cdot \psi_{\Theta^\top} \left(\frac{t_\nu}{2K\varepsilon}\right) \to +\infty.
\end{equation}
Seien   $ C= \frac{\varepsilon}{t_\nu}$ und 
$X = \frac{1}{K \psi_{\Theta^\top} \left(\frac{t_\nu}{2K\varepsilon}\right)}$.
Nach  dem Satz 1 gibt es ein $\pmb{x}\in \mathbb{Z}^m$ mit
$$
||\theta_1 x_1+...+\theta_mx_m - \alpha || < \frac{\varepsilon}{t_\nu}\,\,\,\,\,\text{
und}\,\,\,\,\,
|\pmb{x}|\le X .
$$
Es folgt 
aus (\ref{00}), dass 
 die Ungleichung $ X\le t_\nu$ 
 für jedes hinreichend große $\nu$ gilt. Satz 3 folgt daraus.$\Box$

Es ist 
bekannt, dass Sätze  2  (für $n\ge 2$) und 3 optimal sind. 
Wir kombinieren die Ergebnisse aus den Artikeln von  Jarn\'{\i}k   \cite{J1,J2} und \cite{JE}  und erhalten die folgende Aussage:

     \vskip+0.3cm

{\bf Satz 4}.  

\noindent
(a)
{\it  Sei $ n\ge 2$. Wir betrachten eine  positive Funktion $\varphi (t)$ mit $\lim_{t\to+\infty} \varphi (t) = 0$. Dann gibt es eine
$m\times n$ Matrix
$\Theta$ und ein $\pmb{\alpha}\in \mathbb{R}^n$ mit 
den folgenden Eigenschaften.

\noindent
{\rm 1)}
Das System der Zahlen  $ \theta_{i,j} , 1\le i \le m, \, 1\le j \le n$ ist algebraisch unabhängig.

\noindent
{\rm 2)}  Für die inhomogene Approximationen hat man
\begin{equation}\label{0000}
||{\Theta}\pmb{x} - \pmb{\alpha} || \ge \varphi (|\pmb{x}|),\,\,\,\ \forall \pmb{x}\in \mathbb{Z}^m\setminus \{\pmb{0}\}.
\end{equation}
}

\noindent
(b)
{\it  Seien $ n = 1 $ und $m\ge 2$. 
Wir betrachten eine  positive Funktion $\varphi (t)$ mit $\lim_{t\to+\infty} t\cdot  \varphi (t) = 0$. Dann gibt es eine
$m\times 1$ Matrix
$\Theta = (\theta_1,...,\theta_m)$  und ein $\alpha \in \mathbb{R}$
mit
den folgenden Eigenschaften.

\noindent
{\rm 1)}
Das System der Zahlen  $ \theta_{i} , 1\le i \le m $ ist algebraisch unabhängig.

\noindent
\noindent
{\rm 2)} 
 Die Ungleichung (\ref{0000}) gilt.}

   \vskip+0.5cm
   
{\bf 3. Linear unabhängige Punkte und Fundamentalbereiche.}

   \vskip+0.5cm
   Wir betrachten die Parallelepipede
$$
\Pi (t) =\{ \pmb{z}:\,\, |\Theta\pmb{x} - \pmb{y}|\le t^{-1},\,\,\, |\pmb{x}|\le( \psi_{\Theta^\top} (t))^{-1}\}
,\,\,\,\,
\Pi^*(t) =
\{ \pmb{z}:\,\, |\Theta^\top\pmb{y} - \pmb{x}|\le \psi_{\Theta^\top} (t),\,\,\, |\pmb{y}|\le t\}.
$$
Seien
 $\lambda_1(t)\le...\le\lambda_d(t)$  sukzessive Minima des Parallelepipeds $
\Pi (t) $  und
 $\mu_1(t)\le...\le\mu_d(t)$  sukzessive Minima des Parallelepipeds $
\Pi (t) $. Es ist klar, dass
\begin{equation}\label{odi}
 \mu_1(t)=1.
 \end{equation}
Wir brauchen die Ungleichungen
\begin{equation}\label{MAHL}
{d}^{-1}
\le
\lambda_\nu (t) \mu_{d+1-\nu} (t)
\le (d-1)!,\,\,\,\,\, 1\le \nu \le d,
\end{equation}
die von K. Mahler (siehe \cite{C}, Kapitel V, \S 9) bewiesen wurden.

Dann  gelten die Ungleichungen
$$\lambda_\nu (t) \le \frac{(d-1)!}{\mu_{d+1-\nu}(t)}
\le
\begin{cases}
\frac{(d-1)!}{\mu_d},\,\,\,\,\,\,\,\, \nu=1
\cr
(d-1)!,\,\, \nu \ge 2
\end{cases}
$$  
und
 es gibt  $d$ linear unabhängige primitive\footnote{   Ein Punkt $\pmb{z} =(z_1,...,z_d)\in \mathbb{Z}^d$ heißt primitiv,
 wenn die Gleichung $ {\rm ggT } (z_1,...,z_d) =1$ gilt.}  Gitterpunkte
$$
\pmb{z}(\nu) = (\pmb{x}(\nu), \pmb{y}(\nu))^\top = (x_1(\nu),...,x_m(\nu), y_1(\nu),...,y_n(\nu))^\top \in \mathbb{Z}^d,
\nu  =1 ,...,d 
$$ 
mit
$$
\pmb{z}(1) \in 
{\lambda_{1}(t)} \cdot \Pi (t)
\subset
  \frac{(d-1)!}{\mu_{d}(t)} \cdot \Pi (t);\,\,\,\,\,\,
\pmb{z}(\nu) \in 
{\lambda_{\nu}(t)} \cdot \Pi (t)
\subset
  {(d-1)!}\cdot \Pi (t) , \, 2\le \nu \le d.
$$

Schließlich haben wir die folgenden Ungleichungen bewiesen:
\begin{equation}\label{pupp}
   |\pmb{x}(1)| \le{ \frac{(d-1)!}{\psi_{\Theta^\top}(t)\cdot  \mu_{d}(t)}},\,
   \max_{2\le \nu \le d} |\pmb{x}(\nu)| \le{ \frac{(d-1)!}{\psi_{\Theta^\top}(t)}},\,
     |\overline{\Theta}\,\pmb{z}(1)|\le  {\frac{(d-1)!}{t \cdot \mu_{d}(t)}},
 \,
\max_{2\le \nu \le d}|\overline{\Theta}\,\pmb{z}(\nu)|\le \frac{(d-1)!}{t}.
\end{equation}
    \vskip+0.3cm

{\bf Hilfssatz 1.}
{\it 
Der Parallelepiped $ d!\cdot \Pi (t)$ enthält eine Fundamentalbereich des Gitters $\mathbb{Z}^d$.}

    \vskip+0.3cm

Beweis.
Wir bemerken, dass aus  den Formeln (\ref{odi}) und (\ref{MAHL}) 
$$
\pmb{z}(\nu) \in (d-1)! \cdot \Pi (t)
$$
folgt. Der Parallelepiped
$$
\mathcal{P}=\{ \pmb{z} = l_1\pmb{z}_1+...+ l_d \pmb{z}_d,\,\,\, l_1,...,l_d \in [0,1] \}  
$$
gehört  also zum Parallelepiped $d! \cdot \Pi (t)$.
Es gibt eine
Fundamentalbereich  in $\mathcal{P}$.$\Box$

    \vskip+0.3cm

{\bf Bemerkung 2.} Es gibt Gitterpunkte $\pmb{z}'(1) = \pmb{z}(1), \pmb{z}'(2),...,\pmb{z}'(d)\in d! \cdot \Pi(t)$,
die eine Basis des Gitters $\mathbb{Z}^d$ bilden.

    \vskip+0.3cm
 
{\bf Hilfssatz 2.} {Die Ungleichung 
$$
  (\mu_{d}(t))^{d} \cdot t^n\cdot (\psi_{\Theta^\top}(t))^m \ge (d!)^{-1}
$$
gilt.}
\vskip+0.3cm

Beweis. 
Es gibt $d$ linear unabhängige Gitterunkte $\pmb{z}(\nu) \in \mu_d(t) \cdot \Pi^*(t)$.
Also gilt
$$  2^{d} \, (\mu_{d}(t))^{d} \cdot t^n\cdot   (\psi_{\Theta^\top}(t))^m ={\rm Vol}\, \Pi^*(t) 
\ge \frac{2^{d}}{d!}
$$
und alles ist bewiesen.$\Box$

\vskip+0.3cm
{\bf Hilfssatz 3.}
{\it 
Sei $ T = \frac{1}{\psi_{\Theta^\top}(t)}$.
Dann  gelten die Ungleichungen
\begin{equation}\label{pupp2}
   \max_{1\le \nu \le d} |\pmb{x}(\nu)| \le {(d-1)!}{T},\,\,\,\,\,
     |\overline{\Theta}\,\pmb{z}(1)|\le  \frac{(d-1)! (d!)^{\frac{1}{d}}}{ (Tt)^{\frac{m}{d}}},
 \,\,\,\,\,
\max_{2\le \nu \le d}|\overline{\Theta}\,\pmb{z}(\nu)|\le \frac{(d-1)!}{t}.
\end{equation}
}
\vskip+0.3cm

Beweis. Der Hilfssatz folgt aus (\ref{pupp}). Es bleibt uns nur die zweite Ungleichung  zu beweisen.    
Wir sehen, dass die Ungleichung 
$$
\frac{(\mu_{d}(t))^d\cdot  t^n}{T^m}   \ge
\frac{1}{d!}
$$
wegen des Hilfssatzes 2 gilt.  Nun folgt  die zweite Ungleichung aus (\ref{pupp2})  aus der dritten Ungleichung (\ref{pupp}).$\Box$
\vskip+0.3cm

{\bf Hilfssatz 4.}
{\it
Seien $ n=1$ und $ d = m+1$.
Sei $ T = \frac{1}{\psi_{\Theta^\top}(t)}$.
Nehmen wir an, dass die Ungleichung
\begin{equation}\label{em}
t\cdot \psi_{\Theta^\top}(t)  = \frac{t}{T}\ge M\ge 1
\end{equation}
gilt.

 \vskip+0.3cm
{\rm (a)} Falls  $ \mu_{m+1} (t) \ge W$ 
hat man
\begin{equation}\label{pupp110}
    |\pmb{x}(1)| \le  \frac{m! \,T}{W},
\,\,\,\,\,\,
\max_{2\le \nu\le m+1}
   |\pmb{x}(\nu)| \le  {m! \,T},
\end{equation}
\begin{equation}\label{pupp11}
     |\theta_1x_1(1)+...+\theta_mx_m(1) - y(1)|\le  {\frac{m!}{TW}},
     \,\,\,\,\,
{\max_{2\le \nu \le m+1}}|\theta_1x_1(\nu)+...+\theta_mx_m(\nu) - y(\nu)|\le \frac{m!}{TM},
\end{equation}
und

 \vskip+0.3cm
{\rm (b)} falls $ \mu_{m+1} (t) < W$  hat man

\begin{equation}\label{pupp111}
 \max_{1\le \nu \le m+1}   |\pmb{x}(\nu)| \le  {m! \,T},
\,\,\,\,\,\,
{\max_{1\le \nu \le m+1}}|\theta_1x_1(\nu)+...+\theta_mx_m(\nu) - y(\nu)|\le \frac{m! \, (m+1)!\, W^{m+1}}{T^m}.
\end{equation}
}

\vskip+0.3cm
{\bf Folgerung 1.}
{\it 
Seien  $ n=1, m\ge 2$ und
\begin{equation}\label{WT}
 W =  T^{\frac{m-1}{m+3}} \ge M,
 \,\,\,\,\,
 T_1 =
 \begin{cases}
 T,\,\,\,\,\, \mu_{m+1}(T) \ge W,
 \cr
 TW, \mu_{m+1}(T) < W.
 \end{cases}
 \end{equation}
 Dann   gelten die Ungleichungen
 \begin{equation}\label{pupp110bis}
    |\pmb{x}(1)| \le  \frac{m! \,T_1}{W},
\,\,\,\,\,\,
\max_{2\le \nu\le m+1}
   |\pmb{x}(\nu)| \le  {m! \,T_1},
\end{equation}
\begin{equation}\label{pupp11bis}
     |\theta_1x_1(1)+...+\theta_mx_m(1) - y(1)|\le  {\frac{m! \, (m+1)!}{T_1W}},
     \,
{\max_{2\le \nu \le m+1}}|\theta_1x_1(\nu)+...+\theta_mx_m(\nu) - y(\nu)|\le \frac{m! \, (m+1)!}{T_1M}.
\end{equation}
 }
 
 \vskip+0.3cm
 Im Fall (a) ist die Folgerung klar.
 Wir betrachten  den Fall (b).
 Wir bemerken, dass die Ungleichung
$$
\frac{W^{m+1}}{T^m}= \frac{1}{TW^2}=\frac{1}{T_1W}\le
\frac{1}{T_1M} =
\frac{1}{TWM}
$$
gilt. Die Ungleichungen (\ref{pupp11bis}) folgen daraus im Fall (b).

Beweis von Hilfssatz 4.

Die
Aussage (a) folgt  direkt aus (\ref{pupp}) und (\ref{em}).  Wir beweisen die Aussage (b).
Die erste Ungleichung aus (\ref{pupp111}) ist klar.
Die zweite Ungleichung aus  (\ref{pupp111}) folgt aus dem Hilfssatz 2 und (\ref{pupp}), weil
die Ungleichung
 $
 W^{m+1} t (\psi_{\Theta^\top}(t))^m \ge \frac{1}{(m+1)!}
$
oder
$
 \frac{1}{t} \le \frac{(m+1)! W^{m+1}}{T^m}
$
gilt.$\Box$
 
 \vskip+0.3cm
     
{\bf Bemerkung 3.}

 \vskip+0.3cm
({\bf i}) Aus dem  Hilfssatz 3  folgt, dass es für jedes $\varepsilon > 0$  unendlich viele Matrizen
\begin{equation}\label{Mat}
\mathcal{Z} =
\left(
\begin{array}{ccc}
z_1(1) & ...& z_1(d)\cr
\vdots & \vdots&\vdots\cr
z_d(1) & ...& z_d(d)
\end{array}
\right) =(\pmb{z}(1),...\pmb{z}(d))
\in {\rm SL}_d (\mathbb{Z})
\end{equation}
mit
$$
|\overline{\Theta} \pmb{z} (\nu)| < \varepsilon, \,\,\,
\nu = 1,..., d
$$
gibt.

 \vskip+0.3cm
({\bf ii}) 
Sei $\varphi (t)\downarrow 0$, $ t\to \infty$  eine   fallende Funktion.
In \cite{CheM} wurde gezeigt, dass 
ein Vektor
 $(1,\theta_1,\theta_2) \in \mathbb{R}^3$  
im Fall $ m=1, n=2 , d=3$ 
 mit den folgenden Eigenschaften  existiert.
 
  \vskip+0.3cm
 \noindent
 1)
 Die Zahlen $1, \theta_1,\theta_2$ sind linear unabhängig über $\mathbb{Q}$.
 
  \vskip+0.3cm
 \noindent
 2) 
Für jede Matrix
$$
\left(
\begin{array}{ccc}
x'&x''&x'''\cr
y_{1}'& y_{1}''&y_{1}'''\cr
y_{2}'& y_{2}''&y_{2}'''
\end{array}
\right) \in {\rm SL}_3(\mathbb{Z})
,\,\,\,\,\,  x',x'',x'''\ge 1
$$
gilt die Ungleichung
$$
\max \left\{ \frac{\max_{j=1,2} |x'\theta_j - y_j'|}{\varphi (x')},
 \frac{\max_{j=1,2} |x''\theta_j - y_j''|}{\varphi (x'')},
 \frac{\max_{j=1,2} |x'''\theta_j - y_j'''|}{\varphi (x''')}\right\}
\ge 1.
$$

 \vskip+0.3cm
 
{\bf Bemerkung 4.}

Seien $ n=1, m\ge 2$.  
 
 \vskip+0.3cm
({\bf i})  
Wir nehmen an, dass die Zahlen $1, \theta_1,...,\theta_m$ linear unabhängig über $\mathbb{Q}$ sind.
Sei   $\varepsilon > 0$. Aus der Folgerung 1 und (\ref{J}) folgt, dass  es   unendlich viele Matrizen
 \begin{equation}\label{mat}
\mathcal{Z} = (\pmb{z}(1),...\pmb{z}(d))
\in {\rm SL}_{n+1} (\mathbb{Z})  ,\,\,\,\, \pmb{z}(\nu) = (\pmb{x}(\nu),{y}(\nu))^\top, \,\,\,
\pmb{x}(\nu) \in \mathbb{Z}^m,\,\,\, {y}(\nu)\in \mathbb{Z}
\end{equation}
mit
$$
|\overline{\Theta} \pmb{z} (\nu)|  = ||\theta_1 x_1(\nu)+...+\theta_m x_m(\nu)|| < \frac{\varepsilon}{ |\pmb{x}(\nu)|},\,\,\,
\nu = 1,..., m+1
$$
gibt.

 \vskip+0.3cm
({\bf ii}) 
Wir kennen kein der Bemerkung 3({\bf ii})  ähnliches Ergebnis. Also formulieren wir eine Vermutung:
{\it
für jedes $m\ge 2$  und für jede  Funktion $ \varphi(t)\downarrow 0$  
kann man einen Vektor 
 $\Theta = (\theta_1,...,\theta_m)$  
mit den folgenden Eigenschaften 
konstruieren.\
\ \vskip+0.3cm
 \noindent
 {\rm 1)}
 Die Zahlen $1, \theta_1,..,\theta_m$ sind linear unabhängig über $\mathbb{Q}$.

 \vskip+0.3cm
 \noindent
 {\rm 2)}
Für jede Matrix $\mathcal{Z}$ der Form (\ref{Mat})
mit $ d = m+1$ und $\pmb{z}(\nu) = (x_1(\nu),...,x_m(\nu), y(\nu))^\top\in \mathbb{Z}^{m+1}$  gilt die Ungleichung
$$
\max_{1\le \nu \le m+1}\,\,
\frac{|\pmb{x}(\nu)|\cdot |\theta_1x_1(\nu)+...+\theta_mx_m(\nu)-y(\nu)|}{   \varphi(|\pmb{x}(\nu)|)} \ge 1.
$$
 }

 \vskip+0.3cm

{\bf 4.  Über  primitive Approximationen.}

 \vskip+0.3cm
 
 In 1961 haben  
 J. Chalk und  P. Erd\H{o}s \cite{Ch}
 die folgende Aussage bewiesen.
 
  \vskip+0.3cm

        {\bf Satz 5.}
    {\it Es gibt eine positive Konstante  $ C$  mit der folgenden Eigenschaft.
    Für jedes
    $\theta\in \mathbb{R}\setminus\mathbb{Q} $ und für jedes
       $\alpha$ 
       gibt es unendlich viele primitive Punkte $ (x, y)$, $x\ge 1$ mit
 $$
      |x\theta - \alpha - y|  \le \frac{C}{x}\,\left( \frac{\log x}{\log\log x}
     \right)^2
  .
     $$
    }

      \vskip+0.3cm
    Aus anderer Hand hat der Autor kürzlich \cite{ChaM} das folgende Ergebnis bewiesen.
    
      \vskip+0.3cm
       {\bf Satz 6.} {\it Es gibt  
       eine überabzählbare
        Menge von Paaren 
     $\theta \in \mathbb{R}\setminus \mathbb{Q}$  und $ \alpha\in \mathbb{R}$ mit
    $$
    \inf_{(x,y)\in \mathbb{Z}^2,\,\,   x> 100, \,\,(x,y) =1}  x \,\,\frac{\log \log x}{\sqrt{ \log x}}\,\,|x\theta -\alpha -y|  >0.
    $$}
    
     \vskip+0.3cm
    Zwischen Ober- und Untergrenze klafft  sich hier eine Lücke. Die richtige Ordnung ist unbekannt. Aber im mehrdimensionalen Problem ist alles klar.
       \vskip+0.3cm

{\bf Satz 7}.
{\it

{\rm (a)} Sei $ n \ge 2$.   Sei $ \psi_{\Theta^\top} (t) \neq 0$ für alle $t>1$. Für jedes $\varepsilon>0$ und für jedes  $\pmb{\alpha}$ gibt es unendlich viele primitive Punkte $\pmb{z}\in \mathbb{Z}^d$ mit
$|\overline{\Theta}\pmb{z}-\pmb{\alpha}|<\varepsilon$.

{\rm (b)}  Sei $ n =1 , m\ge 2$. 
Nehmen wir an, dass  die Zahlen $ 1 , \theta_i,\theta_j$  für einigen $ i\neq j$ linear unabhängig  sind.
Für jedes $\varepsilon>0$ und für  jedes ${\alpha}$ gibt es unendlich viele primitive 
Punkte $\pmb{z} = (\pmb{x},y)\in \mathbb{Z}^{m+1}$ mit
$|\theta_1x_1+...+\theta_mx_m-y-{\alpha}|<\frac{\varepsilon}{|\pmb{x}|}$.

 }
 
     \vskip+0.3cm
 Aus dem Satz 4 folgt, dass die Ergebnisse des  Satzes 7 optimal sind.
 
     \vskip+0.3cm

Beweis von Satz 7.

Aus dem Hilfssatz 1 folgt,  dass jede Verschiebung $ \frak{z} + d! \cdot \Pi (t), \frak{z}\in \mathbb{R}^d$ 
des Parallelepipeds $ d! \cdot \Pi (t)$  einen Gitterpunkt enthält. Wir nehmen an, dass $ \overline{\Theta}\,\frak{z}=\pmb{\alpha}$  und $|\frak{z}|\le c$ gelten.
Dann  gibt es für jedes $\pmb{\alpha} \in \mathbb{R}^s$   einen Gitterpunkt $ \pmb{z} = (z_1,...,z_d)$ mit
\begin{equation}\label{ppp}
|\overline{\Theta}\, \pmb{z} - \pmb{\alpha}|
< \frac{d!}{t}, 
\,\,\,\,\,
\text{und}
\,\,\,\,\, |\pmb{z}| \le C(\Theta, \alpha,c) \cdot T.
\end{equation}
Die Vektoren $ \pmb{z}'(1)=\pmb{z}(1), \pmb{z}'(2),...,\pmb{z}'(d)$  bilden eine Basis des Gitters $\mathbb{Z}^d$, darum gibt es ganze Zahlen 
$u_1,...,u_d$ mit
$$
\pmb{z} = u_1\pmb{z}'(1)+...+u_d \pmb{z}'(d).
$$
Es ist klar, dass
für   jedes hinreichend große $c= c(\Theta,\alpha)$
 ein  Gitterpunkt $\frak{z}$ mit
$\max_{2\le j \le d} |u_j|\neq 0$     existiert.
 Nur gilt die Ungleichung
 $$ \max_{1\le j \le d}  |u_j| \le C_1(\Theta, \alpha) T^{d}$$
mit einer positiven Konstante $C_1(\Theta, \alpha) $. Wir nehmen an, dass $u_2\neq 0$.
Betrachten wir die Punkte
$ \pmb{z}+t \pmb{z}(1)$.
Sei $\delta >0$. Dann gibt es ein 
\begin{equation}\label{o}
v =O( u_2^\delta) = O(T^{\delta d})
\end{equation}
mit
$ {\rm ggT} (u_1+v, u_2) = 1$. Der Punkt $ \pmb{z}^* = (\pmb{x}^*, y)^\top= \pmb{z} + v\pmb{z}(1)$ ist primitiv. 

Wir betrachten den Fall (a).
Aus  der zweiten Ungleichung (\ref{pupp2}), (\ref{ppp}) und (\ref{o}) sehen wir, dass
$$
|\overline{\Theta}\, \pmb{z}^* -\pmb{\alpha}|
\le
|\overline{\Theta}\, \pmb{z} -\pmb{\alpha}|+
v |\overline{\Theta}\, \pmb{z}(1)| =  O\left(\frac{1}{t} +\frac{1}{t^{\frac{m}{d}}}
\right)<\varepsilon
$$
für $\delta = \frac{m}{d^2}$  und hinreichend großes $t$ gilt.

Wir betrachten den Fall (b).  Dann
gelten die Ungleichungen (\ref{pupp110bis},\ref{pupp11bis}),
wobei $T_1$ in (\ref{WT}) definiert ist.  
Sei $ \delta = \frac{m-1}{(m+1)(m+3)}$.

Nach (\ref{J}) können wir 
für hinreichend großes $t$     die Ungleichung $ \frac{2}{M} \le\varepsilon$  annehmen. 
Dann
haben wir
$$
|\pmb{x}^*|  = O\left(T_1+ T^{\delta (m+1)}\cdot \frac{T_1}{W}\right) = O(T_1)
$$
und
$$
|\overline{\Theta}\, \pmb{z}^* -\pmb{\alpha}|
\le
|\overline{\Theta}\, \pmb{z} -\pmb{\alpha}|+
v |\overline{\Theta}\, \pmb{z}(1)| =  O\left(\frac{1}{T_1M } + T^{\delta (m+1)} \cdot \frac{1}{T_1W}
\right)<\frac{\varepsilon}{T_1}
$$
für hinreichend großes $t$.

Damit ist alles bewiesen.$\Box$

\end{document}